\documentclass[12pt]{article}
\usepackage[T2A]{fontenc}
\usepackage[cp1251]{inputenc}
\usepackage{amsfonts}
\usepackage{enumerate}
\usepackage{amssymb}
\usepackage{euscript}

\begin{document}

%\begin{titlepage}
%\newpage
%\thispagestyle {empty}

\noindent M.V.Myronyuk

\medskip

\noindent B.Verkin Institute for Low Temperature Physics and
Engineering of the National Academy of Sciences of Ukraine

\medskip

myronyuk@ilt.kharkov.ua

\bigskip

\noindent \textbf{CHARACTERIZATION THEOREMS FOR Q-INDEPENDENT RANDOM
VARIABLES IN BANACH SPACE}

\bigskip

Let $X$ be a Banach space. Following A. Kagan and G. Szekely we
introduce a notion of $Q$-independence for random variables with
values in $X$. We prove analogues of the Skitovich–Darmois and Heyde
theorems for $Q$-independent random variables with values in X.

\bigskip

\bigskip
\bigskip
\bigskip

%\end{titlepage}

%\newpage

\textbf{1. Introduction.} In 1953 V.P.Skitovich and G.Garmois proved
independently that a Gaussian distribution on a real line is
characterized by the independence of two linear forms from $n$
independent random variables (\cite[$\S$3.1]{Kag-Lin-Rao}). In 1962
S.G.Ghurye and I.Olkin generalized the Skitovich-Darmois theorem for
random variables with values in $\mathbb{R}^n$ and nondegenerate
matrices as the coefficients of forms (\cite[$\S$3.2]{Kag-Lin-Rao}).

The Skitovich-Darmois theorem and other characterization theorems
were generalized in different algebraical structures, particularly,
on locally compact Abelian groups and Banach spaces (see e.g.
\cite{Fe7}-\cite{Fe19}).

In 1985 W. Krakowiak generalized the Skitovich–Darmois theorem for
random variables with values in a separable Banach space $X$ and
linear continuous invertible operators as the coefficients of forms
(\cite{Krakowiak}). In the paper \cite{Myr2008} this result was
proved in another way. Besides the Heyde theorem for Banach spaces
was proved in \cite{Myr2008}.

In the paper \cite{KS} A. Kagan and G. Sz\'{e}kely introduced a
notion of $Q$-independence of random variables which generalizes a
notion of independence and proved that some classical theorems hold
true if instead of independence $Q$-independence is considered.
Particularly they proved that the Skitovich-Darmois theorem and the
Cramer theorem about decomposition of a Gaussian distribution hold
true for $Q$-independent random variables. In the paper \cite{Il}
A.Il’inskii studied polynomials which can appear in a notion of
$Q$-independence. In the paper \cite{Rao} B.L.S. Prakasa Rao proved
a generalization of the Kotlarski theorem (see \cite{Kotlarski}) for
$Q$-independent random variables with values in a Hilbert space.
Some generalizations of the Skitovich-Darmois theorem and the Heyde
theorem for $Q$-independent random variables with values in locally
compact Abelian groups were proved in \cite{Fe2017} and
\cite{Myr2019}.

A notion of $Q$-independence can be introduced in a natural way in a
Banach space. We prove that the Skitovich-Darmois theorem and the
Heyde theorem in a Banach space hold true if we consider
$Q$-independence instead of independence (Theorems 1 and 2). These
results generalize the main results of the paper \cite{Myr2008}.
Proofs in these cases follow the scheme of analogous theorems of the
paper \cite{Myr2008}. Also we prove another characterization
theorems (Theorems 3 and 4) for $Q$-independent random variables
with values in a Banach space. In contrast to Theorems 1 and 2,
Theorem 3 and 4 are new not only for $Q$-independent random
variables, but also for independent random variables.

\textbf{2. Notation and definitions.} We use standard facts of the
theory of probability distributions in a Banach space (see, e.g.,
\cite{Vah-Tar-Chob}). Let $X$ be a separable real Banach space and
let $X^*$ be the space dual to it. Denote by ${\rm GL}(X)$ the set
all linear continuous invertible operators in $X$. Denote by $A^*$
the adjoint operator to a linear continuous operator $A$. Denote by
$I$ the identity operator. Denote by $<~x,f>$ the value of a
functional $f\in X^*$ at an element $x\in X$. Denote by
$\mathcal{M}^1(X)$ the convolution semigroup of probability
distributions on $X$. For $\mu \in \mathcal{M}^1(X)$ we denote by
$\overline{\mu}$ the distribution defined by the formula
$\overline{\mu}(E)=\mu(-E)$ for any Borel set $E\subset X$. The
characteristic functional of the probability distribution $\mu$ of a
random variable $\xi$ with values in $X$ is defined by the formula

$$ \widehat\mu(f)=\mathbf{E}[e^{i<\xi,f>}]=\int\limits_X e^{i<x,f>}
d\mu(x), \quad f\in X^*.$$

\noindent Note that
$\widehat{\overline{\mu}}(f)={\overline{\widehat\mu}}(f)$. It is
known that $\widehat\mu(f)$ is positive definite and sequentially
continuous in the pointwise topology (\cite[Ch.4,
\S2]{Vah-Tar-Chob}).

In the paper, unless otherwise specified, we work in a strong
topology, i.e. the topology arising from a norm.

%\medskip

\textbf{Definition 1}. A random variable $\xi$ is called Gaussian if
either it is degenerated or for any $f\in X^*$ the random variable
$<\xi,f>$ is Gaussian.

The following proposition follows immediately from Definition 1.

%\medskip

\textbf{Proposition 1.} \textsl{A distribution $\mu$ on $X$ is
Gaussian iff for any $f\in X^*$ the characteristic function
$\widehat\mu_f(t)$ of a random variable $<\xi,f>$ is a
characteristic function of a Gaussian distribution on $\mathbb{R}$.}

%\medskip

Note that the following proposition is true in a Banach space.

%\medskip

\textbf{Proposition 2} (\cite[Ch. 4, $\S$2]{Vah-Tar-Chob}).
\textsl{Let $X$ be a real Banach space. The characteristic
functional of a Gaussin distribution $\mu\in \mathcal{M}^1(X)$ has
the form
\begin{equation}\label{0.2}
    \widehat\mu(f)= e^{ i<m,f>-{1\over 2}<R f,f>},\ f\in X^*,
\end{equation}
\noindent where $m\in X$, and $R: X^*\longrightarrow X$ is a
symmetric nonnegative operator. The element $m$ is called the mean
value of the distribution $\mu$, and $R$ is called its covariation
operator.}

\textsl{Conversely, if $\mu\in \mathcal{M}^1(X)$ has a
characteristic functional of the form $(\ref {0.2})$, where $m\in X$
is a certain element and $R: X^*\longrightarrow X$ is a certain
symmetric nonnegative operator, then $\mu$ is a Gaussian
distribution in $X$ with mean value $m$ and covariation operator
$R$.}

%\medskip

%%Враховуючи те, що $\widehat\mu(tf)=\widehat\mu_f(t)$ для будь-яких $f\in
%X^*, t\in \mathbb{R}$, отримаємо, що

%\begin{equation}\label{0.1}
%    <m,f>=m_f, \quad <R f,f>=\sigma_f.
%\end{equation}

Let $\psi(f)$ be a function on $X^*$, and let $h$ be an arbitrary
element of $X^*$. We denote by $\Delta_h$ an operator of finite
difference
$$\Delta_h \psi(f)=\psi(f+h)-\psi(f).$$
A function $\psi(f)$ on $X^*$ is called a polynomial if
$$\Delta_{h}^{n+1}\psi(f)=0$$
for some $n$ and for all $f,h \in X^*$.

Let $\xi_1, \dots, \xi_n$ be random variables with values in
separable Banach space $X$. Following A. Kagan and G. Sz\'{e}kely
(\cite{KS}), we say that random variables $\xi_1, \dots, \xi_n$ are
$Q$-independent  if the characteristic functional of the
distribution of the vector $(\xi_1, \dots, \xi_n)$ can be
represented in the form
\begin{equation}\label{i0}
    \widehat\mu_{(\xi_1, \dots, \xi_n)}(f_1, \dots, f_n)={\bf E}[(\xi_1,
f_1)\cdots(\xi_n,
f_n)]=$$$$=\left(\prod_{j=1}^n\widehat\mu_{\xi_j}(f_j)\right)\exp\{q(f_1,
\dots, f_n)\}, \quad f_j\in X^*,
\end{equation}
where $q(f_1, \dots, f_n)$ is a continuous polynomial on $(X^*)^n$
such that $q(0, \dots, 0)=0$.

$Q$-analogue of the Cramer theorem in a Banach space follows
directly from Proposition 1 and the $Q$-analogue of the Cramer
theorem on a real line (see \cite{KS}).

\textbf{Proposition 3.} \textsl{ If a random variable $\xi$ with
values in a real separable Banach space $X$ has a Gaussian
distribution and $\xi=\xi_1+\xi_2$, where $\xi_1, \xi_2$ are
$Q$-independent random variables, then $\xi_1, \xi_2$ have Gaussian
distributions too.}

\medskip

\textbf{3. The Skitovich-Darmois theorem.} We prove an analogue of
the Skitovich-Darmois theorem for $Q$-independent random variables
in a Banach space.

%\medskip

\textbf{Theorem 1}. \textsl{Let $\xi_1,\dots,\xi_n, \ n \ge 2,$ be
$Q$-independent random variables with values in a real separable
Banach space $X$ and with distributions $\mu_j$. Let $A_j, B_j\in
GL(X)$. If linear forms $L_1=A_1\xi_1+\dots+A_n\xi_n$ and
$L_2=B_1\xi_1+\dots+B_n\xi_n$ are $Q$-independent then $\mu_j$ are
Gaussian distributions.}

%\medskip

To prove Theorem 1 we need the following lemmas.

\textbf{Lemma 1}. \textsl{Let $\xi_1,\dots,\xi_n, \ n \ge 2,$ be
$Q$-independent random variables with values in a real separable
Banach space $X$ and with distributions $\mu_j$. Let $A_j, B_j$ be
linear continuous operators in $X$. Linear forms
$L_1=A_1\xi_1+\dots+A_n\xi_n$ and $L_2=B_1\xi_1+\dots+B_n\xi_n$ are
$Q$-independent iff the characteristic functionals
$\widehat\mu_j(y)$ satisfy the
equation
\begin{equation}\label{lem1.1}
    \prod_{j = 1}^{n}{\widehat\mu_j(A_j^* f + B_j^* g)} = \prod_{j
= 1}^{n} {\widehat\mu_j(A_j^*f)} \prod_{j = 1}^{n}
{\widehat\mu_j(B_j^* g)} \exp\{r(f,g)\}, \quad f,g\in X^*,
\end{equation}
\noindent where $r(f,g)$ is a continuous polynomial on $(X^*)^2$,
$r(0,0)=0$.}

The proof of Lemma 1 repeats literally the proof of the similar
lemma for locally compact Abelian groups of the paper \cite{Fe2017}.

\textbf{Lemma 2} (\cite{Myr2008}). \textsl{ If $\mu$ is a Gaussian
distribution in a real separable Banach space $X$ and
$\mu=\mu_1*\mu_2$, then $\mu_1$ and $\mu_2$ are Gaussian
distributions on $X$ too.}

\textbf{Lemma 3} (\cite{Myr2008}). \textsl{Let $X$ be a real
separable Banach space, $\mu\in \mathcal{M}^1(X)$. If
\begin{equation}\label{1.0}
    \widehat\mu(f)=e^{\psi(f)},
\end{equation}
\noindent where $\psi(f)$ is a polynomial, in a certain neighborhood
of zero, then $\mu$ is a Gaussian distribution on $X$.}

\textbf{Proof of Theorem 1.} Putting $\xi'_j=A_j\xi_j$, we reduce
the proof to the case of linear forms $L_1=\xi_1+\dots+\xi_n$ and
$L_2=C_1\xi_1+\dots+C_n\xi_n$, where $C_j\in{\rm GL}(X)$. It follows
from Lemma 1 that the condition of $Q$-independence of $L_1$ and
$L_2$ is equivalent to the statement that the characteristic
functionals $\widehat\mu_j(y)$ satisfy the equation

\begin{equation}\label{5.1}
    \prod_{j = 1}^{n}{\widehat\mu_j(f + C_j^* g)} = \prod_{j
= 1}^{n} {\widehat\mu_j(f)} \prod_{j = 1}^{n} {\widehat\mu_j(C_j^*
g)} \exp\{r(f,g)\}, \quad f,g\in X^*,
\end{equation}

\noindent where $r(f,g)$ is a continuous polynomial on $(X^*)^2$,
$r(0,0)=0$. Thus, the proof of Theorem 1 reduces to the description
of solutions of equation (\ref{5.1}) in the class of characteristic
functionals of probability distributions on $X^*$.

Set $\nu_j=\mu_j*\bar\mu_j$. Then
$\widehat\nu_j(f)=|\widehat\mu_j(f)|^2\geq 0$. It is obvious that
the characteristic functionals $\widehat\nu_j(f)$ satisfy equation
(\ref{5.1}) too. If we prove that $\widehat\nu_j(f)$ are
characteristic functionals of Gaussian distributions, then Lemma 2
implies that $\widehat\mu_j(f)$ are characteristic functionals of
Gaussian distributions too. Therefore we can assume from the
beginning that $\widehat\mu_j(f)\geq 0$.

Since all $\widehat\mu_j(f)\geq 0$, we have $\exp\{r(f,g)\}>0$. The
verification of the fact that $\widehat\mu_j(f)> 0$ for all $f\in
X^*$, $j=1,2,\dots,n$, is the same as in the paper \cite{Myr2008}.

We show that $\widehat\mu_j(f)$ is a characteristic functional of a
Gaussian distribution. Set $\psi_j(f)=-\log \widehat\mu_j(f)$. It
follows from (\ref{5.1}) that

\begin{equation}\label{5.2}
    \sum_{j = 1}^{n}{\psi_j(f + C_j^* g)} = P(f)+Q(g)+r(f,g),
\quad f,g\in X^*,
\end{equation}

\noindent where

$$P(f)=\sum_{j = 1}^{n}{\psi_j(f)},\quad Q(g)=\sum_{j =
1}^{n}{\psi_j(C_j^* g)}.$$

We use the finite differences method to solve (\ref{5.2}) and we
follow the scheme of the proof of the Skitovich-Darmois theorem for
locally compact Abelian groups (see \cite{Fe-SD-2003}). Let $h_n$ be
an arbitrary element of $X^*$. We set $k_n=-(C_n^*)^{-1}h_n$. Then
$h_n+C_n^*k_n=0$. Substitute in (\ref{5.2}) $f+h_n$ for $f$ and
$g+k_n$ for $g$. Subtracting equation (\ref{5.2}) from the resulting
equation we obtain

\begin{equation}\label{5.3}
    \sum_{j = 1}^{n-1} \Delta_{l_{n,j}}{\psi_j(f + C^*_j g)}
    =\Delta_{h_{n}} P(f)+\Delta_{k_{n}} Q(g) +\Delta_{(h_{n},k_{n})}r(f,g),
\quad f,g\in X^*,
\end{equation}

\noindent where $l_{n,j}= h_n+C^*_j k_n=(C^*_j-C^*_n)k_n$,
$j=1,2,\dots,n-1$. We note that the left-hand side of (\ref{5.3})
does not contain the function $\psi_n$. Let $h_{n-1}$ be an
arbitrary element of $X^*$. We set $k_{n-1}=-(C^*_{n-1})^{-1}
h_{n-1}$. Then $h_{n-1}+C^*_{n-1} k_{n-1}=0$. Substitute in
(\ref{5.2}) $f+h_{n-1}$ for $f$ and $g+k_{n-1}$ for $g$. Subtracting
equation (\ref{5.3}) from the resulting equation we obtain

\begin{equation}\label{5.4}
    \sum_{j = 1}^{n-2} \Delta_{l_{n-1,j}} \Delta_{l_{n,j}}{\psi_j(f + C^*_j g)}
    =\Delta_{h_{n-1}} \Delta_{h_{n}} P(f)+\Delta_{k_{n-1}} \Delta_{k_{n}} Q(g)+\Delta_{(h_{n-1},k_{n-1})}\Delta_{(h_{n},k_{n})}r(f,g),
\quad f,g\in X^*,
\end{equation}

\noindent where $l_{n-1,j}= h_{n-1}+C^*_j k_{n-1}=(C^*_j-C^*_{n-1})
k_{n-1}$, $j=1,2,\dots,n-2$. The left-hand side of (\ref{5.4}) does
not contain functions $\psi_n$ and $\psi_{n-1}$. Reasoning
similarly, we get the equation

\begin{equation}\label{5.5}
   \Delta_{l_{2,1}} \Delta_{l_{3,1}}\dots \Delta_{l_{n,1}}{\psi_1(f + C^*_1 g)}=$$
    $$=\Delta_{h_{2}}\Delta_{h_{3}}\dots
    \Delta_{h_{n}} P(f)+\Delta_{k_{2}} \Delta_{k_{3}}\dots \Delta_{k_{n}} Q(g) +\Delta_{(h_{2},k_{2})}\Delta_{(h_{3},k_{3})}\dots\Delta_{(h_{n},k_{n})}r(f,g),
\quad f,g\in X^*,
\end{equation}

\noindent where $h_m$ is an arbitrary element of $X^*$,
$k_m=-(C^*_{m})^{-1} h_m$, $l_{m,1}= h_m+C^*_1
k_m=(C^*_1-C^*_{m})k_m$, $m=2,3,\dots,n$. Let $h_{1}$ be an
arbitrary element of $X^*$. We set $k_{1}=-(C^*_{1})^{-1} h_{1}$.
Then $h_{1}+C^*_{1} k_{1}=0$. Substitute in (\ref{5.5}) $f+h_{1}$
for $f$ and $g+k_{1}$ for $g$. Subtracting equation (\ref{5.5}) from
the resulting equation we obtain

\begin{equation}\label{5.6}
   \Delta_{h_{1}}\Delta_{h_{2}}\dots
    \Delta_{h_{n}} P(f)+\Delta_{k_{1}} \Delta_{k_{2}}\dots \Delta_{k_{n}} Q(g) +\Delta_{(h_{1},k_{1})}\Delta_{(h_{2},k_{2})}\dots\Delta_{(h_{n},k_{n})}r(f,g)=0,
\quad f,g\in X^*.
\end{equation}

\noindent Let $h$ be an arbitrary element of $X^*$. Substitute in
(\ref{5.6}) $f+h$ for $f$. Subtracting equation (\ref{5.6}) from the
resulting equation we obtain

\begin{equation}\label{5.7}
   \Delta_{h}\Delta_{h_{1}}\Delta_{h_{2}}\dots
    \Delta_{h_{n}} P(f)+\Delta_{(h,0)}\Delta_{(h_{1},k_{1})}\Delta_{(h_{2},k_{2})}\dots\Delta_{(h_{n},k_{n})}r(f,g)=0,
\quad f, g\in X^*.
\end{equation}

Since $r(f,g)$ is a polynomial, we have
\begin{equation}\label{5.7.1}
  \Delta^{l+1}_{(h, k)}r(f,g)=0,
\quad f, g\in X^*,
\end{equation}
for some $l$ and arbitrary $h$ and $k$ from $X^*$.

We apply $\Delta^{l+1}_{(h, k)}$ to (\ref{5.7}). Taking into account
 (\ref{5.7.1}), we get

\begin{equation}\label{5.7.2}
   \Delta^{l+2}_{h}\Delta_{h_{1}}\Delta_{h_{2}}\dots
    \Delta_{h_{n}} P(f)=0,
\quad f\in X^*.
\end{equation}

Note that $f,h$ and $h_m$, $m=1,..,n$, are arbitrary elements of
$X^*$. We can put in (\ref{5.7}) $h_1=\dots=h_n=h$. Then

\begin{equation}\label{5.8}
   \Delta_{h}^{n+l+2}P(f)=0, \quad f,h\in X^*.
\end{equation}

\noindent Thus, $P(f)$ is a polynomial on $X^*$. We denote
$\gamma=\mu_1*\dots*\mu_n$. Then
 $\widehat\gamma(f)=\prod_{j = 1}^{n}\widehat\mu_j(f)$. Thus,
$\widehat\gamma(f)=e^{-P(f)}$. It follows from Lemma 3 that
$\widehat\gamma(f)$ is a characteristic functional of a Gaussian
distribution. Then it follows from Proposition 3 that
$\widehat\mu_j(f)$ is a characteristic functional of a Gaussian
distribution. Theorem 1 is proved. $\blacksquare$

\medskip

\textbf{4. The Heyde theorem.} We prove an analogue of the Heyde
theorem for $Q$-independent random variables in a Banach space.

%\medskip

\textbf{Theorem 2}. \textsl{Let $\xi_1,\dots,\xi_n, \ n \ge 2,$ be
$Q$-independent random variables with values in a real separable
Banach space $X$ and with distributions $\mu_j$. Let $A_j, B_j
\in{\rm GL}(X)$ such that $B_i A_i^{-1} \pm B_j A_j^{-1} \in{\rm
GL}(X)$ for all $i \ne j$. If the conditional distribution of the
linear form $L_2=B_1\xi_1+\dots+B_n\xi_n$ given
$L_1=A_1\xi_1+\dots+A_n\xi_n$ is symmetric, then $\mu_j$ are
Gaussian distributions.}

%\medskip

To prove Theorem 2 we need the following lemma.

\textbf{Lemma 4}. \textsl{Let $\xi_1,\dots,\xi_n, \ n \ge 2,$ be
$Q$-independent random variables with values in a real separable
Banach space $X$ and with distributions $\mu_j$. Let $A_j, B_j
\in{\rm GL}(X)$. The conditional distribution of the linear form
$L_2=B_1\xi_1+\dots+B_n\xi_n$ given $L_1=A_1\xi_1+\dots+A_n\xi_n$ is
symmetric iff the characteristic functional $\widehat\mu_j(y)$
satisfy equation
\begin{equation}\label{lem4.1}
    \prod_{j = 1}^{n}{\widehat\mu_j(A^*_jf + B^*_j g)} =
    \prod_{j = 1}^{n}{\widehat\mu_j(A^*_jf - B^*_j g)} \exp\{r(f,g)\} ,
\quad f,g\in X^*,
\end{equation}
\noindent where $r(f,g)$ is a continuous polynomial on $(X^*)^2$,
$r(0,0)=0$.}

The proof of Lemma 4 repeats literally the proof of the similar
lemma for locally compact Abelian groups of the paper \cite{Fe2017}.

\textbf{Proof of Theorem 2.} As in the proof of Theorem 1, we show
that the proof of Theorem 2 reduces to the case of
$L_1=\xi_1+\dots+\xi_n$ and $L_2=C_1\xi_1+\dots+C_n\xi_n$, where
$C_j\in{\rm GL}(X)$ such that $C_i \pm C_j \in{\rm GL}(X)$ for all
$i \ne j$, and characteristic functional $\widehat\mu_j(f)\geq 0$.
It follows from Lemma 4 that the condition of symmetry of the
conditional distribution of $L_2$ given $L_1$ is equivalent to the
statement that the characteristic functionals $\widehat\mu_j(y)$
satisfy the equation

\begin{equation}\label{6.1}
    \prod_{j = 1}^{n}{\widehat\mu_j(f + C^*_j g)} =
    \prod_{j = 1}^{n}{\widehat\mu_j(f - C^*_j g)} \exp\{r(f,g)\} ,
\quad f,g\in X^*,
\end{equation}

\noindent where $r(f,g)$ is a continuous polynomial on $(X^*)^2$,
$r(0,0)=0$.

Since $\widehat\mu_j(0)=1$ and $\widehat\mu_j(f)$ are sequentially
continuous, there exists a neighborhood $U$ of zero such that
$\widehat\mu_j(y)>0$ for all $y\in U$. We choose in $U$ a symmetric
neighborhood $V$ of zero such that
$$
\sum_{j=1}^{4n} \lambda_j (V) \subset U,
$$
where $\lambda_j \in \{I, C^*_1, \dots,  C^*_n\}$.

We take the logarithm of equation (\ref{6.1}) in the neighborhood
$V$ and get
\begin{equation}\label{6.2}
\sum_{j = 1}^{n} \psi_j(f + C^*_j g)  = \sum_{j = 1}^{n} \psi_j(f -
C^*_j g) + r(f,g), \quad f, g \in V,
\end{equation}
where $\psi_j(f)=-\log \widehat\mu_j(f)$. We use the finite
differences method to solve (\ref{6.2}) and we follow the scheme of
the proof of the classical Heyde (see
\cite[$\S$13.4.1]{Kag-Lin-Rao}).

Let $h_n$ be an arbitrary element of $V$. Substitute in (\ref{6.2})
$f+C_n^* h_n$ for $f$ and $g+h_n$ for $g$. Subtracting equation
(\ref{6.2}) from the resulting equation we obtain

\begin{equation}\label{6.3}
\sum_{j = 1}^{n} \Delta_{l_{n,j}} \psi_j(f + C^*_j g)  = \sum_{j =
1}^{n-1} \Delta_{m_{n,j}} \psi_j(f - C^*_j g) +
\Delta_{(h_{n},h_{n})} r(f,g), \quad f, g \in V,
\end{equation}

\noindent where $l_{n,j}=(C_n^*+C^*_j) h_n$, $m_{n,j}=(C_n^*-C^*_j)
h_n$. Note that the right-hand side of (\ref{6.3}) does not contain
the function $\psi_n$. Let $h_{n-1}$ be an arbitrary element of $V$.
Substitute in (\ref{6.3}) $f+C^*_{n-1}h_{n-1}$ for $f$ and
$g+h_{n-1}$ for $g$. Subtracting equation (\ref{6.3}) from the
resulting equation we obtain

\begin{equation}\label{6.4}
\sum_{j = 1}^{n} \Delta_{l_{n-1,j}}\Delta_{l_{n,j}} \psi_j(f + C^*_j
g)  = \sum_{j = 1}^{n-2} \Delta_{m_{n-1,j}}\Delta_{m_{n,j}} \psi_j(f
- C^*_j g)+ \Delta_{(h_{n-1},h_{n-1})}\Delta_{(h_{n},h_{n})} r(f,g),
\quad f, g \in V,
\end{equation}

\noindent where $l_{n-1,j}=(C_{n-1}^*+C^*_j) h_{n-1}$,
$m_{n,j}=(C_{n-1}^*-C^*_j) h_{n-1}$. The right-hand side of
(\ref{6.4}) does not contain the functions $\psi_n$ and
$\psi_{n-1}$. Reasoning similarly, we get the equation

\begin{equation}\label{6.5}
\sum_{j = 1}^{n}
\Delta_{l_{1,j}}\Delta_{l_{2,j}}\dots\Delta_{l_{n,j}} \psi_j(f +
C^*_j g)  =
\Delta_{(h_{1},h_{1})}\Delta_{(h_2,h_2)}\dots\Delta_{(h_{n},h_{n})}
r(f,g) , \quad f, g \in V,
\end{equation}

\noindent where $h_i$ is an arbitrary element of $V$,
$l_{i,j}=(C_{i}^*+C^*_j) h_{i}$, $i=1,2,\dots,n$.

Let $k_n$ be an arbitrary element of $V$. Substitute in (\ref{6.5})
$f+C_n^* k_n$ for $f$ and $g-k_n$ for $g$. Subtracting equation
(\ref{6.5}) from the resulting equation we obtain

\begin{equation}\label{6.6}
\sum_{j = 1}^{n-1} \Delta_{b_{n,j}}
\Delta_{l_{1,j}}\Delta_{l_{2,j}}\dots\Delta_{l_{n,j}} \psi_j(f +
C^*_j g)  = \Delta_{(C_n^*
h_n,-h_n)}\Delta_{(h_{1},h_{1})}\Delta_{(h_2,h_2)}\dots\Delta_{(h_{n},h_{n})}
r(f,g), \quad f, g \in V,
\end{equation}

\noindent where $b_{n,j}=(C_n^*-C^*_j) k_n$.  Note that the
left-hand side of (\ref{6.6}) does not contain the function
$\psi_n$. Let $k_{n-1}$ be an arbitrary element of $V$. Substitute
in (\ref{6.6}) $f+C^*_{n-1}k_{n-1}$ for $f$ and $g-k_{n-1}$ for $g$.
Subtracting equation (\ref{6.6}) from the resulting equation we
obtain

\begin{equation}\label{6.7}
\sum_{j = 1}^{n-2} \Delta_{b_{n-1,j}}\Delta_{b_{n,j}}
\Delta_{l_{1,j}}\Delta_{l_{2,j}}\dots\Delta_{l_{n,j}} \psi_j(f +
C^*_j g)  $$$$= \Delta_{(C^*_{n-1}k_{n-1},-k_{n-1})}\Delta_{(C_n^*
h_n,-h_n)}\Delta_{(h_{1},h_{1})}\Delta_{(h_2,h_2)}\dots\Delta_{(h_{n},h_{n})}
r(f,g), \quad f, g \in V,
\end{equation}

\noindent where $b_{n-1,j}=(C_{n-1}^*-C^*_j) k_{n-1}$. The left-hand
side of (\ref{6.7}) does not contain the functions $\psi_n$ and
$\psi_{n-1}$. Reasoning similarly, we get the equation

\begin{equation}\label{6.8}
\Delta_{b_{2,1}}\Delta_{b_{3,1}}\dots\Delta_{b_{n,1}}
\Delta_{l_{1,1}}\Delta_{l_{2,1}}\dots\Delta_{l_{n,1}} \psi_1(f +
C^*_1 g)  $$$$= \Delta_{(C^*_{2}k_{2},-k_{2})}\dots\Delta_{(C_n^*
h_n,-h_n)}\Delta_{(h_{1},h_{1})}\dots\Delta_{(h_{n},h_{n})} r(f,g),
\quad f, g \in V,
\end{equation}

\noindent where $k_i$ is an arbitrary element of $V$,
$b_{i,1}=(C_{i}^*-C^*_1) k_{i}$, $i=2,3,\dots,n$.

Since $r(f,g)$ is a polynomial, the equality (\ref{5.7.1}) holds
true for some $l$ and arbitrary $h$ and $k$ from $X^*$. Taking it
into account, we apply the operator $\Delta^{l+1}_{(h, k)}$ to the
both sides of (\ref{6.8}). Putting $g=0$ in the obtained equation,
we get

\begin{equation}\label{6.8.1}
\Delta^{l+1}_{h}\Delta_{b_{2,1}}\Delta_{b_{3,1}}\dots\Delta_{b_{n,1}}
\Delta_{l_{1,1}}\Delta_{l_{2,1}}\dots\Delta_{l_{n,1}} \psi_1(f)  =
0, \quad f \in V,
\end{equation}

Since $k_i, h_i$ are arbitrary elements of $V$, $C_i^* \pm C_j^*
\in{\rm GL}(X^*)$ for all $i \ne j$, and $l_{i,j}=(C_{i}^*+C^*_j)
h_{i}$, $b_{i,1}=(C_{i}^*-C^*_1) k_{i}$, we get in some neighborhood
$W$ of zero that

\begin{equation}\label{6.9}
\Delta_{h}^{2n+l} \psi_1(f)  = 0, \quad f, h \in W.
\end{equation}

\noindent We obtain that $\widehat\mu_1(f)=e^{-\psi_1(f)}$, where
$\psi_1(f)$ is a polynomial, in the neighborhood $W$ of zero. It
follows from Lemma 3 that $\mu_1$ is a Gaussian distribution.
Similarly we obtain that all $\mu_j$ are Gaussian distributions.
$\blacksquare$

If in Theorem 2 $n=2$, then the condition on coefficients of linear
forms can be relaxed. Indeed, it is shown in the proof of Theorem 2
that we can assume that for $n=2$ the conditional distribution of
$L_2=C_1\xi_1+C_2\xi_2$ given $L_1=\xi_1+\xi_2$ is symmetric. Apply
the operator $C_1^{-1}$ to $L_2$. It is obvious that the conditional
distribution of $C_1^{-1} L_2$ given $L_1$ is symmetric too. Thus we
can assume that $L_1=\xi_1+\xi_2$ and $L_2=\xi_1+C\xi_2$.

The following statement holds true.

\textbf{Theorem 3}. \textsl{Let $\xi_1,\xi_2$ be $Q$-independent
random variables with values in a real separable Banach space $X$
and with distributions $\mu_j$. Let $C \in{\rm GL}(X)$ such that
$C+I \in{\rm GL}(X)$. If the conditional distribution of the linear
form $L_2=\xi_1+C\xi_2$ given $L_1=\xi_1+\xi_2$ is symmetric, then
$\mu_j$ are Gaussian distributions.}

To prove Theorem 3 we need the following lemma.

\textbf{Lemma 5}. \textsl{ Let $\xi_1,\xi_2$ be $Q$-independent
random variables with values in a real separable Banach space $X$
and with distributions $\mu_j$. Let $C_j \in{\rm GL}(X)$. If the
conditional distribution of the linear form $L_2=C_1\xi_1+C_2\xi_2$
given $L_1=\xi_1+\xi_2$ is symmetric, then
$L'_1=(C_1+C_2)\xi_1+2C_2\xi_2$ and $L'_2=2C_1\xi_1+(C_1+C_2)\xi_2$
are $Q$-independent.}

\textbf{Proof.} By Lemma 4 the characteristic functionals
$\widehat\mu_j(y)$ satisfy the equation

\begin{equation}\label{l4.1}
    \widehat\mu_1(f+C_1^* g)\widehat\mu_2(f+C_2^* g)
    =\widehat\mu_1(f-C_1^* g)\widehat\mu_2(f-C_2^* g) \exp
    \{r(f,g)\},\quad f,g\in X^*,
\end{equation}
where $r(f,g)$ is a continuous polynomial on $(X^*)^2$, $r(0,0)=0$.

Putting $f=C_2^* h, g=-h$ and then $f=-C_1^* h, g=h$ in equation
(\ref{l4.1}), we get

\begin{equation}\label{l4.2}
    \widehat\mu_1((C_2^* - C_1^*)h)=
    \widehat\mu_1((C_1^*+C_2^*)h)\widehat\mu_2(2C_2^*h) \exp
    \{r(C_2^* h,-h)\},\quad h\in X^*,
\end{equation}

\begin{equation}\label{l4.3}
    \widehat\mu_2((C_2^*-C_1^*)h)=
    \widehat\mu_1(-2C_1^*h)\widehat\mu_2(-(C_1^*+C_2^*)h) \exp
    \{r(-C_1^* h,h)\},\quad h\in X^*.
\end{equation}

Let $k, l\in X^*$. Put $f=C_1^* l+C_2^* k$, $g=k+l$ in (\ref{l4.1}).
We obtain

\begin{equation}\label{l4.5}
    \widehat\mu_1((C_1^*+C_2^*)k+2C_1^* l)
    \widehat\mu_2(2C_2^* k+(C_1^*+C_2^*)l)
    =$$$$=\widehat\mu_1((C_2^*-C_1^*)k)
    \widehat\mu_2(-(C_2^*-C_1^*)l) \exp
    \{r(C_1^* l+C_2^* k, k+l)\},\quad k,l\in X^*.
\end{equation}

Taking into account (\ref{l4.2}) and (\ref{l4.3}), we can rewrite
equation (\ref{l4.5}) in the form

\begin{equation}\label{l4.6}
    \widehat\mu_1((C_1^*+C_2^*)k+2C_1^* l)
    \widehat\mu_2(2C_2^* k+(C_1^*+C_2^*)l)
    =$$$$=
    \widehat\mu_1((C_1^*+C_2^*)k)\widehat\mu_2(2C_2^*k) \widehat\mu_1(2C_1^*l)\widehat\mu_2((C_1^*+C_2^*)l)
    \times $$$$ \times \exp
    \{r(C_1^* l+C_2^* k, k+l)+r(C_2^* k,-k)+r(-C_1^* l,l)\},\quad k,l\in X^*.
\end{equation}

It follows from Lemma 1 and equation (\ref{l4.6}) that the linear
forms $L'_1=(C_1+C_2)\xi_1+C_2\xi_2$ and
$L'_2=2C_1\xi_1+(C_1+C_2)\xi_2$ are $Q$-independent. $\blacksquare$

\textbf{Proof of Theorem 3.} It follows from Lemma 5 that the linear
forms $L'_1=(I+C)\xi_1+2C\xi_2$ and $L'_2=2\xi_1+(I+C)\xi_2$ are
$Q$-independent. Since $I+C \in{\rm GL}(X)$, all coefficients of
linear forms $L'_1$ and $L'_2$ are linear continuous invertible
operators. It follows from Theorem 1 that Theorem 3 holds true.
$\Box$

\textbf{Remark 1.} Suppose that in Theorem 3 a Banach space $X$ is
reflexive. It is not difficult to verify that the condition $I+C
\in{\rm GL}(X)$ can be relaxed to the condition $Ker (I+C)=\{0\}$,
and the assertion of the theorem does not change. Indeed, in this
case $\overline{Im (I+C^*)}=X^*$. It follows from the proof of
Theorem 1 that only this condition is essential.

Since the independence of random variables implies their
$Q$-independence, it follows from Theorem 3 the following statement.

\textbf{Corollary 1.} \textsl{Let $\xi_1,\xi_2$ be independent
random variables with values in a real separable Banach space $X$
and with distributions $\mu_j$. Let $C \in{\rm GL}(X)$ such that
$C+I \in{\rm GL}(X)$. If the conditional distribution of the linear
form $L_2=\xi_1+C\xi_2$ given $L_1=\xi_1+\xi_2$ is symmetric, then
$\mu_j$ are Gaussian distributions.}

\medskip

\textbf{5. The independence of the sample mean and the residue
vector.} It is well known that the Gaussian distribution on a real
line is characterized by the independence of the sample mean and the
residue vector. Sufficiency of this statement is obvious. Necessity
follows from the Geary-Lucacs-Laga theorem
(\cite[$\S$4.2]{Kag-Lin-Rao}). This characterization of a Gaussian
distribution was generalized on locally compact Abelian groups in
the paper \cite{FeMyr-2010}. In this section we prove an analogue of
this characterization for $Q$-independent random variables with
values in a Banach space.

%\medskip

\textbf{Theorem 4}. \textsl{Let $\xi_1,\dots,\xi_n, \ n \ge 2,$ be
$Q$-independent identically distributed random variables with values
in a real separable Banach space $X$ and with a distribution $\mu$.
Put ${S={1\over n}\sum_{j=1}^n\xi_j}$ and ${\bf V}=(\xi_1-S, \dots,
\xi_n-S)$. If ${S}$ and ${\bf V}$ are $Q$-independent then $\mu$ is
a Gaussian distribution.}

%\medskip

To prove Theorem 4 we need the following lemma.

\textbf{Lemma 6}. \textsl{Let $\xi_1,\dots,\xi_n, \ n \ge 2,$ be
$Q$-independent identically distributed random variables with values
in a real separable Banach space $X$ and with a distribution $\mu$.
Put ${S={1\over n}\sum_{j=1}^n\xi_j}$ and ${\bf V}=(\xi_1-S, \dots,
\xi_n-S)$. If ${S}$ and ${\bf V}$ are $Q$-independent then the
characteristic functional $\widehat\mu(f)$ satisfy equation
\begin{equation}\label{eq1}
    \prod\limits_{j=1}^{n}\widehat\mu({f\over n}+ g_j -{1\over n}(g_1+\dots+g_n))$$$$=\widehat\mu^n({f\over n})
\prod\limits_{j=1}^{n}\widehat\mu( g_j -{1\over n}(g_1+\dots+g_n))
\exp\{r(f,g_1,\dots,g_n)\}, \quad f,g,g_j, \in X^*,
\end{equation}
\noindent where $r(f,g_1,\dots,g_n)$ is a continuous polynomial on
$(X^*)^2$, $r(0,0,\dots,0)=0$.}

\textbf{Proof.} Note that $S$ and ${\bf V}$ are $Q$-independent iff
the equality
\begin{equation}\label{l1.2}
 \mathbf{E}[(S, f)(\mathbf{V}, (g_1,\dots,g_n))]=
 \mathbf{E}[(S, f)]
 \mathbf{E}[(\mathbf{V}, (g_1,\dots,g_n))] \exp\{q_1(f,g_1,\dots,g_n)\}
\end{equation}
holds for all $u, v_1,\dots,v_n \in X^*$, where
$q_1(f,g_1,\dots,g_n)$ is a continuous polynomial on $(X^*)^2$,
$q_1(0,0,\dots,0)=0$. Since random variables $\xi_1,\dots,\xi_n$ are
$Q$-independent, the left-hand side of (\ref{l1.2}) can be
represented in the form

 $$\mathbf{E}[(S, f)(\mathbf{V}, (g_1,\dots,g_n))]=$$$$
\mathbf{E}[({1\over n}(\xi_1+\dots+\xi_n), f)((\xi_1-{1\over
n}(\xi_1+\dots+\xi_n), \dots,\xi_n-{1\over n}(\xi_1+\dots+\xi_n)),
(g_1,\dots,g_n))]=$$
$$ \mathbf{E}[\prod_{j=1}^{n}(\xi_j, {f\over n}+g_j{{-{1\over n}(g_1+\dots+g_n)}})=
\prod_{i=1}^{n} \widehat\mu({f\over n}+ g_i{{ -{1\over
n}(g_1+\dots+g_n)}}) \exp\{q_2(f,g_1,\dots,g_n)\}. $$

Analogously we transform the right-hand side of (\ref{l1.2}):

$$ \mathbf{E}[(S, f)]
 \mathbf{E}[(\mathbf{V}, (g_1,\dots,g_n))]=$$$$=
 \mathbf{E}[({1\over n}(\xi_1+\dots+\xi_n),{{f}})]
 \mathbf{E}[((\xi_1-{1\over n}(\xi_1+\dots+\xi_n),\dots,\xi_n-{1\over n}(\xi_1+\dots+\xi_n)),
 (g_1,\dots,g_n))]=$$
$$=\prod_{j=1}^{n}\mathbf{E}[(\xi_j, {f\over n})] \exp\{q_3(f,g_1,\dots,g_n)\}
 \mathbf{E}[\prod_{j=1}^{n}(\xi_j, g_j-{1\over n}(g_1+\dots+g_n))=$$$$=
 \prod_{j=1}^{n}\mathbf{E}[(\xi_j, {f\over n})]\exp\{q_3(f,g_1,\dots,g_n)\} \prod_{j=1}^{n}\mathbf{E}[(\xi_j,
 g_j-{1\over n}(g_1+\dots+g_n))] \exp\{q_4(f,g_1,\dots,g_n)\}=$$$$=
\widehat\mu^n ({f\over n}) \prod_{i=1}^{n} \widehat\mu( g_i-{1\over
n}(g_1+\dots+g_n))
\exp\{q_3(f,g_1,\dots,g_n)+q_4(f,g_1,\dots,g_n)\}.$$

Put

$$\exp\{r(f,g_1,\dots,g_n)\}= \exp\{q_1(f,g_1,\dots,g_n)+q_3(f,g_1,\dots,g_n)+q_4(f,g_1,\dots,g_n)-q_2(f,g_1,\dots,g_n)+\}. $$

Lemma is proved.
 $\Box$

\textbf{Proof of Theorem 4.} It follows from Lemma 6 that the
function $\widehat\mu(f)$ satisfies equation (\ref{eq1}).
Substituting $f=nh, g_1=g, \ g_2=-g, g_3=\dots=g_n=0$ in
(\ref{eq1}), we obtain
\begin{equation}\label{eq24}
\widehat\mu(h+g)\widehat\mu(h-g)\widehat\mu^{n-2}(h)=
\widehat\mu^n(h)|\widehat\mu(g)|^2 \exp\{r(nh,g,-g,0,\dots,0)\} ,
\quad h,g\in X^*,
\end{equation}

\noindent where $r(f,g_1,\dots,g_n)$ is a continuous polynomial on
$(X^*)^2$, $r(0,0,\dots,0)=0$.

Let $n=2$. Equation (\ref{eq24}) takes the form of equation
(\ref{5.1}) for $n=2$ and $C_1=I$, $C_2=-I$. Since it follows from
Lemma 1 that the linear forms $L_1=\xi_1+\xi_2$ and
$L_2=\xi_1-\xi_2$ are independent, Theorem 4 follows in this case
from the proof of Theorem 1.

Let $n\ge 3$. We can write equation (\ref{eq24}) in the form
\begin{equation}\label{eq25}
\widehat\mu(f+g)\widehat\mu(f-g)\widehat\mu^{n-2}(f)=
\widehat\mu^n(f)|\widehat\mu(g)|^2 \exp\{\widetilde{r}(f,g)\}, \quad
f,g\in X^*,
\end{equation}
where $\widetilde{r}(f,g)$ is a continuous polynomial on $(X^*)^2$,
$r(0,0)=0$.

It follows from (\ref{eq25}) that the set
$$
B=\{f\in Y:\widehat\mu(f)\ne 0\}
$$
is an open subgroup of $X^*$. Thus $B=X^*$ and $\widehat\mu(f)\ne 0$
on $X^*$. As in the proof of Theorem 1, we can show that we can
assume from the beginning that $\widehat\mu_j(f)> 0$. Then we also
have that $r(f,g_1,\dots,g_n)>0$.

Put $\psi(f)=-\log\widehat\mu(f)$. Taking the logarithm of
(\ref{eq1}) and multiplying all arguments of functions in
(\ref{eq1}) on $n$, we get
\begin{equation}\label{eq6}
    \sum\limits_{j=1}^{n}\psi(f+n g_j -(g_1+\dots+g_n))$$$$=n\psi(f)+
\sum\limits_{j=1}^{n}\psi(n g_j -(g_1+\dots+g_n))
+r(f,g_1,\dots,g_n), \quad g,g_j, \in X^*.
\end{equation}

We use the finite difference method to solve equation (\ref{eq6}).
Let $h_1$ be an arbitrary element of $X^*$. Substitute $f+h_1$ for
$f$ and $g_j+h_1$ for $g_j$, $j=1, 2, \dots, n$ in equation
(\ref{eq6}). Subtracting equation (\ref{eq6}) from the resulting
equation we obtain
\begin{equation}\label{eq7}
    \sum\limits_{j=1}^{n}\Delta_{h_1}\psi(f+n g_j -(g_1+\dots+g_n))=
    n\Delta_{h_1}\psi(f) + \Delta_{(h_1,h_1,\dots,h_1)} r(f,g_1,\dots,g_n),
\quad f,h_1, g_j, \in X^*.
\end{equation}
Нехай $h_2$ --- це довільний елемент $X^*$. Покладемо $f+h_2$
замість $f$ та $g_1+h_2$ замість $g_1$ в рівнянні (\ref{eq6}).
Віднімаючи (\ref{eq7}) з отриманого рівняння, ми маємо
\begin{equation}\label{eq8}
    \Delta_{nh_2}\Delta_{h_1}\psi(f+n g_1 -(g_1+\dots+g_n))=$$$$=
    n\Delta_{h_2}\Delta_{h_1}\psi(f)+ \Delta_{(h_2,h_2,0,\dots,0)}\Delta_{(h_1,h_1,\dots,h_1)} r(f,g_1,\dots,g_n),
\quad f,h_1, h_2, g_j, \in X^*.
\end{equation}
Let $h_3$ be an arbitrary element of $X^*$. Substitute $f+h_3$ for
$f$ and $g_2+h_3$ for $g_2$ in equation (\ref{eq6}). Subtracting
equation (\ref{eq8}) from the resulting equation we obtain
\begin{equation}\label{eq9}
   n\Delta_{h_3}\Delta_{h_2}\Delta_{h_1}\psi(f)+$$$$+ \Delta_{(h_3,0,h_3,0,\dots,0)}\Delta_{(h_2,h_2,0,\dots,0)}\Delta_{(h_1,h_1,\dots,h_1)} r(f,g_1,\dots,g_n) =0.
\quad f,h_1, h_2, h_3, g_j, \in X^*.
\end{equation}

Since $r(f,g_1,\dots,g_n)$ is a polynomial, we have
\begin{equation}\label{eq9.1}
  \Delta^{l+1}_{(k, k_1,\dots,k_n)}r(f,g_1,\dots,g_n)=0,
\quad f, g_j\in X^*,
\end{equation}
for some $l$ and arbitrary elements $k$ and $k_j$ from $X^*$.

Applying the operator $\Delta^{l+1}_{(k, k_1,\dots,k_n)}$ to
(\ref{eq9}) and taking into account (\ref{eq9.1}), we get

\begin{equation}\label{eq9.2}
   \Delta^{l+1}_{k}\Delta_{h_3}\Delta_{h_2}\Delta_{h_1}\psi(f) =0.
\quad f,h_1, h_2, h_3, k \in X^*.
\end{equation}

Put $h_1=h_2=h_3=k=h$ in (\ref{eq9.2}). We obtain
\begin{equation}\label{eq9a}
   \Delta_{h}^{l+4}\psi(f)=0,
\quad f, h \in X^*,
\end{equation}
i.e. $\psi(y)$ is a continuous polynomial. Thus,
$\widehat\mu(f)=e^{-\psi(f)}$. It follows from Lemma 3 that
$\widehat\mu(f)$ is a characteristic functional of a Gaussian
distribution. Theorem 4 is proved. $\blacksquare$

Since the independence of random variables implies their
$Q$-independence, it follows from Theorem 4 the following statement.

\textbf{Corollary 2.} \textsl{ Let $\xi_1,\dots,\xi_n, \ n \ge 2,$
be independent identically distributed random variables with values
in a real separable Banach space $X$ and with a distribution $\mu$.
Put ${S={1\over n}\sum_{j=1}^n\xi_j}$ and ${\bf V}=(\xi_1-S, \dots,
\xi_n-S)$. If ${S}$ and ${\bf V}$ are independent then $\mu$ is a
Gaussian distribution.}

\medskip

%\newpage

\end{document}